\documentclass[12pt]{article}
\usepackage{amsmath}
\usepackage{amssymb}
\usepackage[mathscr]{eucal}
\usepackage[usenames]{color}
\usepackage[utf8]{inputenc}
\usepackage{url}
\usepackage{longtable}
\usepackage{graphicx}
\usepackage{epstopdf}
\usepackage{subfigure}

\usepackage{authblk}

\newcommand{\ws}{\om^*}
\newcommand{\ls}{\lambda^*}
\newcommand{\bs}{b^*}
\newcommand{\ps}{\pi^*}
\newcommand{\om}{\omega}

\usepackage{url}
\usepackage[
bookmarks=true,
backref=true,
colorlinks=true,
linkcolor=red,
citecolor=blue,
urlcolor=green]{hyperref}

\newtheorem{theorem}{Theorem}
\newtheorem{lemma}{Lemma}
\newtheorem{proposition}{Proposition}

\newtheorem{corollary}{Corollary}

\numberwithin{equation}{section}

\def\om{\omega}
\definecolor{darkolivegreen}{rgb}{0.333333, 0.419608, 0.1843140}

\def\om{\omega}

\title{\Large Delayed Keen Model with Inflation
}

\author[]{Ali Tolga Dincer\thanks{dincer15@itu.edu.tr}}
\author[]{Sevgi Harman\thanks{harman@itu.edu.tr}}
\author[]{Seyma Gonul\thanks{gonul20@itu.edu.tr}}
\author[]{Ayse Tiryakioglu\thanks{aysetiryakioglu@itu.edu.tr}}
\author[]{Cihangir Ozemir\thanks{Corresponding author, ozemir@itu.edu.tr}}
\affil[]{Department of Mathematics, Istanbul Technical University, Istanbul, Türkiye}

\begin{document}

\date{}

\maketitle

\begin{abstract}
Keen's   model describes the dynamics between wage share, employment rate and debt ratio. In literature, the model was extended to represent the effects of inflation and also the speculative money flow. Based on the inflationary model, we take into account a time delay in the inflation term which stands for the period before the effects of inflation are seen. We observe that, the delayed system may experience a Hopf bifurcation and exhibit cyclic behavior around an equilibrium point, although the non-delayed model is stable under the same conditions.

\vspace{0.5cm}
\end{abstract}

\section{Introduction}
Keen model \cite{Keen1995,Keen2013} is one of the well-known models in macroeconomics. It is considered in \cite{Grasselli2015}  as the three dimensional system 
\begin{equation}
	\label{okm}
	\begin{aligned}
		\dot{\om} &= \om\big[\Phi(\lambda)-\alpha\big], \\
		\dot{\lambda} &= \lambda\big [g(\pi)-\alpha-\beta\big], \\
		\dot{b} &= \kappa(\pi) - \pi - b. g(\pi)
	\end{aligned}
\end{equation}
where $\om$ is the wage share, $\lambda$ is the employment rate and $b$ is the debt ratio. $\Phi$ is called a Phillips curve which determines the relation between wage rate and employment rate and $\kappa$ is the investment function which determines the relation between debt and profit rate $\pi$. $g$ is growth function  and $\alpha$, $\beta$ are constants. Various improvements have been done on this model, such as the introduction of an inflation function $Z(w)$ according to a price dynamics given in \cite{Keen2013} and \cite{desai1973growth}  by Grasselli and Huu \cite{Grasselli2015},  which is expressed as
\begin{equation}
	\label{kmwi}
	\begin{aligned}
		\dot{\om} &= \om\big[\Phi(\lambda)-\alpha-(1-\gamma)Z(\om)\big], \\
		\dot{\lambda} &= \lambda\big [g(\pi)-\alpha-\beta\big], \\
		\dot{b} &= \kappa(\pi) - \pi - b\big[Z(\om)+g(\pi)\big].
	\end{aligned}
\end{equation}
We want to enhance the approach of Grasselli and Huu, by assuming that the inflation function   $Z(\om(t))$ affects  the dynamics with  a certain delay $\tau\geq 0$, which reflects the time  before the effects of inflation are seen. Therefore we propose the delayed version of the Keen model with inflation \eqref{kmwi} as 
\begin{equation}
	\label{main}
	\begin{aligned}
		\dot{\om} &=\om \big[\Phi(\lambda) - \alpha - (1-\gamma) Z(\om(t-\tau))\big], \\
		\dot{\lambda} &= \lambda \big[g(\pi) - \alpha - \beta\big], \\
		\dot{b} &= \kappa(\pi) - \pi - b \big[Z(\om(t-\tau))+g(\pi)\big].
	\end{aligned}
\end{equation}
We perform stability analysis at a fixed point and show that the system experiences a Hopf bifurcation. 

A comprehensive analysis of the Keen model \eqref{okm} is carried out in \cite{Grasselli2012}. They show that the system \eqref{okm} has two economically relevant equilibrium points. One is economically desirable with finite debt, strictly positive wages and employment rate. The other equilibrium point is with zero wage share and employment ratio, and an infinite debt ratio, which is therefore an undesirable one. They give stability conditions of these equilibrium points. Further, they extend the model by a modification on the Ponzi  variable which was suggested by \cite{keen2009policy}.

The analysis in  \cite{Grasselli2015} reveals that,  unlike the original Keen model, the Keen  model with inflation \eqref{kmwi} exhibits three different equilibrium points: A \emph{"good"} one with finite and nonzero values, and two \emph{"bad"} equilibrium points, one with infinite debt share, one with zero employment rate \cite{Grasselli2015}. The article further improves \eqref{kmwi} to a four-component dynamics by also considering  a speculative money variable.

As is well  known, in many fields, such as biology \cite{Zhang2014,Faria2001,X.ZhangandH.Zhao2017}, chemistry \cite{Beta2003,J.WeinerF.W.Schneider1989} and engineering \cite{Agrawal1997,Agrawal2000,Lu2018}, delay is commonly encountered. Inevitably, economic mechanisms will also depend on delays. For example, changes in the money supply will not cause immediate changes in the economy, there will always be a period of lag. The production cycle has long and short phases. There is always a delay in price changes. Hence, the delay differential equations will model a realistic economic mathematical model \cite{DeCesare2005,Fanti2007}.

Ma and  Chen \cite{jun2001study} work on a financial model of which variables are interest rate, investment demand and price index.  The analyses of  \cite{Gao2009} and \cite{Ding2012} investigate the case when there is a delay in the investment demand in this system. 
Zhang \cite{Zhang2019}  builds its model based on the Chen's system and adds delay to the price index.  The studies \cite{Chen2008} and \cite{Son2011}
consider a more general case, in which all the  above mentioned  variables also affect the dynamics with different delays and feedback strengths. For those cases with delay, the financial models experience a Hopf bifurcation at the equilibrium  points for a certain range of the parameters.

In this paper, the dynamics of the macroeconomic model \eqref{main} will be analysed, in which  a delay occurs in the inflation term. This dynamics has its roots in the Keen model  \eqref{okm} first proposed in \cite{Keen1995,Keen2013}. The Keen model takes three main variables into account: wages, employment and private debt, inspired by the Goodwin Business Cycle Model \cite{Goodwin1967}, which describes the macroeconomic fluctuations of a single closed commodity economy.
In \cite{Grasselli2012} and \cite{Grasselli2015}, price dynamics are added to the model. The new system remains three-dimensional, as price level enters the system implicitly through wage share and debt ratio. The stability condition of the equilibrium point that already exists in the original Keen model, must be modified in order to express it in nominal terms and include inflationary or deflationary systems. The explanation is actually the same, but as the price level dynamically increases, certain conditions are weakened, while others are strengthened  \cite{Grasselli2015,Grasselli2012}.

Delayed forms of the Keen model and its  variants have not been analyzed so far in the existing literature, to the best of the our knowledge. By the assumption that the effect of the inflation contributes to the dynamics with a time delay, we propose system \eqref{main} and work on its stability. We find that under suitable conditions on the parameters, although the non-delayed model is stable at an equilibrium point, the delayed system  undergoes a Hopf bifurcation at a certain value of the delay. Occurrence of a Hopf bifurcation represents fluctuation in the dynamic macroeconomic variable(s), which may correspond to real world observation, hence may provide a more realistic model.

We present the analysis and the results in the following section and we conclude by comments and discussions. We also include an Appendix section, in which we perform the calculations for the characterization of the Hopf bifurcation.

\section{Analysis}

\subsection{Information on the model}

For the sake of completeness and to provide an understanding of the evolution in \eqref{main}, we mention the explanation of the terms the system is containing, originally made in \cite{Grasselli2015}.
In Keen model \cite{Keen2013}, the output $Y$ and capital $C$ have a constant capital-to-output ratio $\nu$
\begin{equation}\label{LitEq1}
	Y=\frac{C}{\nu}
\end{equation}
and capital is defined to change dynamically as
\begin{equation}\label{LitEq2}
	\dot{C} =   I-\delta C
\end{equation}
where $ I$ is the real investment by firms and $\delta$ is the constant depreciation rate. Ref. \cite{Grasselli2015} says that firms obtain funds both internally and externally and internal funds are the net profits, that is,
\begin{equation}\label{LitEq3}
	\Pi = pY - W - rB,
\end{equation}
where $p$ is the price level for goods, $W$ is the wage bill, $r$ is the interest rate and $B$ is the net borrowing (debt) of firms from banks. The change of the borrowing that firms take from banks depends on the nominal investment as
\begin{equation}\label{LitEq4}
	\dot{B} = p  I - \Pi.
\end{equation}
According to \cite{Grasselli2015}, Keen \cite{Keen1995} defines investment as $ I=\kappa(\pi)Y$ where $\kappa(\pi)$ is a function of the net profit share satisfying the conditions in  \cite{Grasselli2012}
\begin{equation}
	\begin{aligned}
		&	\kappa'(\pi)>0 \quad  \text{on}   \quad (-\infty,\infty), \\
		&	\lim_{\pi \to -\infty} \kappa(\pi) = \kappa_0 < \nu(\alpha+\beta+\delta) <  \lim_{\pi \to +\infty} \kappa(\pi),\\
		&		\lim_{\pi \to -\infty} \pi^2 \kappa'(\pi) = 0
	\end{aligned}
\end{equation} 
and the profit rate is
\begin{equation}\label{LitEq5}
	\pi = \frac{\Pi}{pY} = 1-\om-rb
\end{equation}
with
\begin{equation}\label{LitEq6}
	\om = \frac{W}{pY}, 
	\qquad  \quad 
	b = \frac{B}{pY}
\end{equation}
denoting the wage and firm debt shares, respectively. Let us mention that,  following \cite{Grasselli2012} and \cite{Grasselli2015} here we will take $\kappa(\pi)$ as
\begin{equation}
	\kappa(\pi) = \kappa_0 + e^{\kappa_1 + \kappa_2 \pi}
\end{equation}
where $\kappa_0, \kappa_1$ and $\kappa_2$ are suitable constants.

Let the total workforce be  $N$, the number of employed workers be  $l$. The productivity of a worker  $a$, the employment rate as $\lambda$, and the nominal wage rate  $ w$ are defined in  \cite{Grasselli2015} as follows
\begin{equation}
	a = \frac{Y}{l},
	\; \; \; \;
	\lambda = \frac{l}{N} = \frac{Y}{aN}
	\; \; \; \;
	w = \frac{W}{l}.
\end{equation}
Dynamic changes of productivity and workforce are assumed to obey exponential growth, that is,
\begin{equation}
	\frac{\dot{a}}{a} = \alpha, \qquad 
	\frac{\dot{N}}{N} = \beta
\end{equation}
where  $\alpha$ and $\beta$ are constants, leading to the employment rate dynamics
\begin{equation}
	\frac{\dot{\lambda}}{\lambda} = \frac{\dot{Y}}{Y} - \alpha-\beta.
\end{equation}
In Keen model with inflation \cite{Grasselli2015}, the dynamics of wage rate is related to the employment rate by a classic Phillips curve as
\begin{equation}
	\begin{aligned}
		& \frac{\dot{ w}}{ w} = \Phi(\lambda) + \gamma Z,\\
		& Z(\om)=\frac{\dot p}{p}=\eta_p (\xi \om -1)
	\end{aligned}
\end{equation}
for constants $0\leq \gamma \leq 1$, $\eta_p >0$, $\xi \geq 1$.  The function $\Phi$ is a continuously differentiable function defined on $[0, 1)$ satisfying
\begin{equation}
	\begin{aligned}
		&	\Phi'(\lambda)>0 \quad  \text{on}   \quad (0,1), \\
		&	\Phi(0)<\alpha,\\
		&		\lim_{\lambda \to 1^{-}} \Phi(\lambda) =\infty. 
	\end{aligned}
\end{equation}
Let us note that we shall consider it to be 
\begin{equation}
	\Phi(\lambda) = \frac{\Phi_1}{(1 - \lambda)^2} - \Phi_0
\end{equation}
as in \cite{Grasselli2015}.
Finally, growth rate $g$ is defined as
\begin{equation}
	g(\pi) = \frac{\dot{Y}}{Y} = \frac{\kappa(\pi)}{\nu}-\delta.
\end{equation}
The above definitions lead to the dynamics \eqref{okm} if $\dot p \equiv 0$ and to \eqref{kmwi} if $\dot p \not\equiv 0$.   For a more detailed discussion, we refer to \cite{Keen1995,Keen2013} and  \cite{Grasselli2015,Grasselli2012}. 

\subsection{Delayed Keen model with inflation}

We introduce a time delay to the system in the inflation function. Our suggestion is that the change of inflation in an economy will not show its effect directly, instead the effect will be observed after a certain amount of time. Let us recall that the system we are working on is \eqref{main} where $Z(\om(t))=\eta_p (\xi \om(t)-1)$ and $\pi=1-\om-rb$. Let $(\ws,\ls,\bs)$ be an equilibrium point of this system. As $\ws=\om^*(t)=\om^*(t-\tau)$, the equilibrium points are found from the solution of the system 
\begin{subequations}
	\begin{eqnarray}
		\ws .\big[\Phi(\ls) - \alpha - (1-\gamma) Z(\ws)\big]=0,&   \label{eqa}\\
		\ls .\big[g(\ps) - \alpha - \beta\big]=0,&  \label{eqb}\\
		\kappa(\ps) - \ps -\bs .\big[Z(\ws)+g(\ps)\big]=0,& \label{eqc}\\
		\ps+\ws+r\bs =1.&  \label{eqd}
	\end{eqnarray}
\end{subequations}
If $\ws=\ls=0$,  we find the equlibrium point $E_1(0,0,b_1^*)$ where $b_1^*$ is a solution to the equation
\begin{equation}
	\kappa(1-rb_1^*)-(1-rb_1^*)=b_1^*[g(1-rb_1^*)-\eta_p]
\end{equation}
If $\ws\neq 0$, $\ls =0$, we obtain the equlibrium point   $E_2(\om_2^*,0,b_2^*)$ where 
\begin{equation}
	\om_2^*=\frac{\Phi(0)-\alpha}{(1-\gamma) \eta_p \xi}+\frac{1}{\xi}
\end{equation}
and  $b_2^*$ is a root of the equation 
\begin{equation}
	\kappa(1-\om_2^*-rb_2^*)-(1-\om_2^*-rb_2^*)=b_2^*[g(1-\om_2^*-rb_2^*)+\eta_p (\xi \om_2^*-1)].
\end{equation}
In case $\ws=0$, $\lambda \neq 0$, searching for the equilibrium point $E_3(0,\ls_3,\bs_3)$, \eqref{eqa} is satisfied. \eqref{eqb} requires $\ps_3 = g^{-1}(\alpha+\beta)$ and \eqref{eqd} gives $\bs_3=(1-\ps_3)/r$. If this value of $\bs_3$ also satisfies \eqref{eqc}, that is, if 
\begin{equation}
	\kappa(\ps_3) - \ps_3=\bs_3 .\big[g(\ps_3)-\eta_p\big],\\
\end{equation}
then $E_3(0,\ls_3,\bs_3)$ is an equilibrium point with arbitrary  $\ls_3$. 

Let $(\ws,\ls)\neq (0,0)$. Then necessarily 
\begin{equation}
	\label{ProfitEquilibriumPoint}
	\pi^* = g^{-1}(\alpha+\beta) = \kappa^{-1}(\nu(\alpha+\beta+\delta))
\end{equation}
and we find $E_4(\ws,\ls,\bs)$ 	with
\begin{subequations}
	\label{OldEquilibriumPoints}
	\begin{eqnarray}
		\om^* &=& 1-\pi^*-rb^*, \label{eq2a}\\
		\lambda^* &=& \Phi^{-1} (\alpha+(1-\gamma)\  Z(\om^*)), \label{eq2b} \\
		b^* &=& \frac{\kappa(\pi^*)-\pi^*}{Z(w^*)+\alpha+\beta}\, \label{eq2c}
	\end{eqnarray}
\end{subequations}
as the equilibrium points. As explained  in \cite{Grasselli2015}, \eqref{eq2a} and \eqref{eq2c} together give that $\ws$ is a solution to the equation  
\begin{equation}
	a_0(\ws)^2+a_1 \ws +a_2=0
\end{equation}
with 
\begin{equation}
	\begin{aligned}
		a_0 &= \xi\eta_p, \\
		a_1 &= \alpha+\beta- \eta_p-\xi \eta_p(1-\ps),\\
		a_2 &= (\eta_p-\alpha-\beta)(1-\ps)+r(\kappa(\pi^*)-\pi^*).
	\end{aligned}
\end{equation}

We shall assume that at the equilibrium point $E=E_4(\ws,\ls,\bs)$  we have  $\ws, \ls\in(0,1)$ and illustrate the stability behavior and occurrence of the Hopf bifurcation at only this equilibrium point.
The Jacobian matrices of the system \eqref{main} in the non-delayed and delayed parameters are respectively

\begin{equation}
	\label{Jacobian}
	J_0=
	\begin{bmatrix}
		\Phi(\lambda)-\alpha-(1-\gamma)\eta_p (\xi \om_\tau-1) & \om\Phi'(\lambda) & 0 \\
		-\lambda g'(\pi) & g(\pi)-\alpha-\beta & -r\lambda g'(\pi) \\
		1-\kappa'(\pi)+bg'(\pi) & 0 & r(1-\kappa'(\pi)+bg'(\pi))-(Z(\om_{\tau})+g(\pi))
	\end{bmatrix},
\end{equation}

\begin{equation}
	J_{\tau}=
	\begin{bmatrix}
		(\gamma - 1) \om \eta_p \xi &\quad & 0 &\quad& 0 \\
		0 && 0 && 0 \\
		-b \eta_p \xi && 0 && 0
	\end{bmatrix}
\end{equation}
where $\om_\tau=\om(t-\tau)$.
At $E$ we have $g(\ps) = \alpha+\beta$. The Jacobians become
\begin{equation}
	J_0\bigr\vert_{E}=
	\begin{bmatrix}
		0 & K_1 & 0 \\
		-K_2 & 0 & -rK_2 \\
		K_3  & 0 & K_4
	\end{bmatrix},
	\quad
	J_{\tau}\bigr\vert_{E}=
	\begin{bmatrix}
		K_0 &\quad  0 &\quad  0 \\
		0 &\quad  0 &\quad  0 \\
		-\bs \eta_p \xi	 & \quad  0 &\quad  0
	\end{bmatrix}
\end{equation}
with the terms
\begin{equation}
	\begin{aligned}
		K_0 &= (\gamma - 1) \ws \eta_p \xi, \\
		K_1 &= \ws \Phi'(\ls),\\
		K_2 &= \ls \frac{\kappa'(\ps)}{\nu},\\
		K_3 &= 1-\kappa'(\ps)+\bs g'(\ps),\\
		K_4 &= rK_3-(\alpha+\beta+Z(\ws)). \\			
	\end{aligned}
\end{equation}
The eigenvalue equation is 
\begin{align}
	P(x)&=\vert J_0+e^{-x\tau}J_{\tau}-xI \vert\bigr\vert_{E}   \nonumber \\
	    &=
	\begin{vmatrix}
		e^{-x\tau} K_{0} - x &\quad  K_1 & \quad  0 \\
		-K_2 &\quad -x &\quad -rK_2 \\
		K_3 - e^{-x\tau}  \bs \eta_p \xi &\quad 0 &\quad K_4 - x
	\end{vmatrix}=0,    \label{char}
\end{align}
where $x$ represents the eigenvalue.
If $\tau=0$, the characteristic equation  is
\begin{equation}\label{ndeq}
	P_0(x)=-x^3 +(K_0+K_4)x^2 -(K_0 K_4+K_1 K_2)x-K_1 K_2 K_5=0
\end{equation}
where $K_5=\alpha+\beta+Z(\ws)-\bs r\eta_p\xi$. The eigenvalue equation $P_0(x)=0$ has all of its roots on the left half-plane if the Routh-Hurwitz criterion is satified.
\begin{lemma}
	\label{RHLemma}
	If the conditions listed as
	\begin{equation}\label{RH}
		\begin{aligned}
			K_0 + K_4 < 0, \\
			K_1 K_2 K_5 > 0,\\
			K_1 K_2 K_5 + (K_0 + K_4) (K_0 K_4 + K_1 K_2) < 0
		\end{aligned}
	\end{equation}
	hold, then the system fulfills the Routh-Hurwitz criterion. This means that provided these criteria are held, all three roots of the characteristic equation \eqref{ndeq} have negative real parts and thus the equilibrium point is stable.
\end{lemma}
There is another condition in \cite{Grasselli2015} as $ K_0 K_4+K_1 K_2  > 0$, which actually follows from the above three. In general, characteristic equation \eqref{char} takes the form
\begin{equation}  \label{deq}
	P(x) = -x^3 + K_4 x^2 - K_1 K_2 x - K_1 K_2 K_7   
	     + e^{-x \tau} [ K_0 x^2- K_0 K_4  x-r K_1 K_2 K_6  ]=0 	
\end{equation}
with $K_6 = -\bs \eta_p \xi$, $K_7=\alpha+\beta+Z(\ws)$.
Equation \eqref{deq} reduces to equation \eqref{ndeq} when $\tau=0$. 

Let $x=i\mu$ be a solution of \eqref{deq}. It  can be separated to real and imaginary parts as
\begin{align}
	\label{Real and Imaginary Eqs}
\mu^3-K_1 K_2 \mu &= 	K_0 K_4 \mu \cos(\mu\tau) - (rK_1 K_2 K_6  + K_0 \mu^2) \sin(\mu\tau)   \\
  -\mu^2 K_4-K_1 K_2K_7  &=	(rK_1 K_2 K_6  + K_0 \mu^2)\cos( \mu\tau)  + K_0 K_4 \mu \sin(\mu\tau) 	. 
\end{align}
Here an elimination to
\begin{equation}\label{eqw}	
	\mu^6 + p\mu^4 + q\mu^2 + \tilde{r} = 0
\end{equation}
is possible where
\begin{align}
		p &=  K_{4}^{2} - K_0^2-2 K_1 K_2,   \nonumber\\
		q &= K_1^2 K_2^2- K_0^2 K_4^2  + 2K_1 K_2 K_4 K_7  - 2r K_0 K_1 K_2 K_6,  \nonumber\\
		\tilde{r} &= K_1^2 K_2^2 (K_7^2- r^2K_6^2).
\end{align}
By the substitution $z=\mu^2$, equation \eqref{eqw} becomes
\begin{equation}\label{eqh}
	h(z)=z^3 + pz^2 + qz + \tilde{r} =0.
\end{equation}
The question now is to determine whether this cubic equation has any roots $z=\mu^2>0$ or not. We use the results available in \cite{Song2004}. Clearly, if $\tilde r <0$, then \eqref{eqh} has at least one positive root. In case $\tilde r \geq 0$, then the derivative
\begin{equation}\label{eqht}
	h'(z)=3z^2 + 2pz + q
\end{equation}
has the zeros  $\displaystyle z_1^*=\frac{1}{3}(-p+\sqrt{\Delta})$ and $\displaystyle z_2^*=\frac{1}{3}(-p-\sqrt{\Delta})$ where $\Delta=p^2-3q$. Without giving the details, we express the result available in \cite{Song2004} as follows.
\begin{proposition}
	The following statements are true for the roots of Eq. \eqref{eqh}. 
	
	\renewcommand{\labelenumi}{\roman{enumi}}
	\label{Positive root conditions}
	\begin{enumerate}
		\item If $\tilde{r}<0$ then there exists at least one positive root.
		\item If $\tilde{r}\geq0$ and $\Delta=p^2-3q\leq 0$, then there is no positive root.
		\item If $\tilde{r}\geq0$ and $\Delta=p^2-3q>0$, then there exists a positive root if and only if $z_1^*>0$ and $h(z_1^*)\leq 0$.
	\end{enumerate}
\end{proposition}
Let us assume that \eqref{eqh} has positive roots. Without loss of generality, we can assume there are three positive roots $z_1$, $z_2$, $z_3$. So \eqref{eqw} has the roots $\mu_1=\sqrt{z_1}$, $\mu_2=\sqrt{z_2}$, $\mu_3=\sqrt{z_3}$.

From equations \eqref{Real and Imaginary Eqs}, we have

\begin{equation}\label{eqtau}
	\tau_k^j=\frac{1}{\mu_k}\left\{\cos^{-1}\left(-\frac{\mu^2K_0 K_1 K_2 K_4 + r \mu^2 K_1 K_2 K_4 K_6 +  \mu^2K_0 K_1 K_2 K_7 +r K_1^2 K_2^2 K_6 K_7}{\mu^4 K_0^2+\mu^2 K_0^2 K_4^2+2r\mu^2 K_0 K_1 K_2 K_6+r^2 K_1^2 K_2^2 K_6^2}\right)+2j\pi\right\}
\end{equation}

\noindent where $k=1,2,3$; $j=0,1,2,\ldots$; then, $\mp i \mu_k$ is a pair of pure imaginary roots to \eqref{deq} with $\tau_k^j$. Define
\begin{equation}
	\tau_0=\tau_{k_0}^0=\min_{k\in \{1,2,3\}}\{\tau_k^0\}, \quad \mu_0 = \mu_{k_0}, \quad z_0 =\mu_0^2.
\end{equation}
\begin{lemma}\label{ruan}
	Consider the exponential polynomial equation
	\begin{align}
		P(x,e^{-x \tau_{1}},\ldots, e^{-x \tau_{m}})&= 
		\, x^{n}+p_{1}^{(0)} x^{n-1}+\ldots+p_{n-1}^{(0)} x+p_{n}^{(0)}   \nonumber\\
		&+[p_{1}^{(1)} x^{n-1}+\ldots +p_{n-1}^{(1)} x +p_{n}^{(1)}]e^{-x \tau_{1}}    \\
		&+\ldots    \nonumber\\
		&+[p_{1}^{(m)} x^{n-1}+\ldots+p_{n-1}^{(m)} x +p_{n}^{(m)}]e^{-x \tau_{m}}=0,   \nonumber
	\end{align}
	where $\tau_{i} \geq 0$ $(i=1,2,...,m)$ and $p_{j}^{(i)}$ $(i=0,1,2,...,m; j=1,2,...,n)$ are constants. As $(\tau_{1},\tau_{2},...,\tau_{m})$ vary, the sum of the order of the zeros of $P(x, e^{-x \tau_{1}},...,e^{-x \tau_{m}})$ on the open right half-plane can change only if a zero appears on or crosses the imaginary axis.
\end{lemma}
\begin{corollary}
	Suppose the conditions \eqref{RH} are satisfied.
	\begin{itemize}
		\item[(i)] If $\tilde{r}\geq0$ and $\Delta=p^2-3q\leq 0$, then, the roots of \eqref{deq} are with negative real parts for all $\tau\geq 0$. 
		\item[(ii)] If (a) $\tilde{r}<0$ or (b) $\tilde{r}\geq0$,   $\Delta=p^2-3q>0$, $z_1^*>0$ and $h(z_1^*)\leq 0$, then, the roots of \eqref{deq} are with negative real parts for all $\tau\in [0,\tau_0)$. 
	\end{itemize}
\end{corollary}
\begin{lemma}
	Assume that $h'(z_{0})\neq0$. Then, $\mp i\mu_{0}$ is a simple (i.e., not multiple) pure imaginary root of the equation \eqref{deq} when $\tau=\tau_{0}$. Additionally, the  transversality condition 
	\begin{equation}
		\frac{d(\mathrm{Re}(x(\tau)))}{d\tau}\Big\vert_{\tau=\tau_{0}}\neq0
	\end{equation}
	is satisfied	and the sign of $d(\mathrm{Re}(x(\tau)))/d\tau\vert_{\tau=\tau_{0}}$ is consistent with that of $h'(z_{0})$ if the conditions (ii) of Corollary 1  are satisfied.
\end{lemma}
\noindent
\textit{Proof}:  	Define $R(x), Q(x)$ as
\begin{equation}
	\begin{aligned}
		R(x) &= -x^3 + K_4 x^2 - K_1 K_2 x - K_1 K_2 K_7 ,\\
		Q(x)&= K_0 x^2 - K_0 K_4 x - r K_1 K_2 K_6. 
	\end{aligned}
\end{equation}
Then equation \eqref{deq} becomes
\begin{equation}\label{RQ}
	R(x) + Q(x) e^{-x\tau} = 0.
\end{equation}
Let us evaluate \eqref{RQ} and its complex conjugate at $x=i\mu$ to get
\begin{equation}
	\begin{aligned}
		&R(i\mu) + Q(i\mu) e^{-i\mu\tau} = 0,\\
		&\bar R(i\mu) + \bar Q(i\mu) e^{i\mu\tau} = 0.\\ 
	\end{aligned}
\end{equation}
If we eliminate the exponential term between these two equations, what we get is nothing but \eqref{deq}
\begin{equation}\label{hw2}
	h(\mu^2) = R(i\mu) \bar{R}(i\mu) - Q(i\mu) \bar{Q}(i\mu)=0.
\end{equation}
The rest of the proof is by working on \eqref{hw2} through standard calculations, at the end of which one obtains
\begin{align}
	&\frac{d(\mathrm{Re}(x(\tau)))}{d\tau}\Big\vert_{\tau=\tau_{0},x=i\mu_{0}}  \nonumber \\
	&\quad =\frac{\mu_{0}^{2}h'(\mu_{0}^{2})}{\vert R'(i\mu_0)e^{i\mu_0 \tau_0}+ Q'(i\mu_0)-\tau_0 Q(i\mu_0)\vert^{2}}\neq 0.  
\end{align}
The details of the calculations are available  in  Lemma 7 of \cite{Zhang2019} and, was used in Lemma 9 of \cite{Demirci2020} for the same purpose. Furthermore, if the conditions (ii) of Corollary 1 are satisfied, then  we also have,
$\displaystyle \frac{d(\mathrm{Re}(x(\tau)))}{d\tau}\Big\vert_{\tau=\tau_{0}}> 0 $, as otherwise for $\tau <\tau_0$, $\mathrm{Re}(x(\tau))$ would be decreasing and  assume positive values.

Evaluation of the Lyapunov coefficient and the related constants characterizing the Hopf bifurcation are presented in the Appendix section. Now we can give our theorem.
\begin{theorem}
	Let $E=E_4(\ws,\ls,\bs)$ be the equilibrium point for which the Routh-Hurwitz conditions in Lemma 1 are satisfied  and let $\tau_0$, $z_0$ be defined as before. 
	\begin{itemize}
		\item  If $\tilde{r}\geq0$ and $\Delta=p^2-3q\leq 0$, then $E$ is asymptotically stable for all $\tau \geq 0$.		
		\item If (a) $\tilde{r}<0$ or (b) $\tilde{r}\geq0$,   $\Delta=p^2-3q>0$, $z_1^*>0$,  $h(z_1^*)\leq 0$ and $h'(z_0)\neq 0$, then   $E$ is an asymptotically stable equilibrium point when $\tau \in [0, \tau_0 )$. $E$ is unstable for $\tau \in (\tau_0, \tau_1 )$,  where $\tau_1$ is the first root of \eqref{deq} after $\tau_0$. The system \eqref{main} undergoes a Hopf bifurcation at  $E$ when $\tau = \tau_0$.
	\end{itemize}
\end{theorem}

\subsection{Numerical simulations}

In this part of the article we shall illustrate the above results by a concrete example.
Following \cite{Grasselli2015}, we take the Philips curve $\Phi$ as 
\begin{equation}
	\Phi(\lambda) = \frac{\Phi_1}{(1 - \lambda)^2} - \Phi_0
\end{equation}
with 
\begin{equation}
	\Phi_0=0.04340277,  \quad  \Phi_1 = 0.00006944.
\end{equation}
The other constants are chosen to be
\begin{equation}
	\alpha=0.025, \quad  \beta=0.02, \quad \delta=0.01 \quad  \nu=3, \quad  r=0.03.
\end{equation}
The constants of the  function $\kappa(\pi)$  are
\begin{equation}
	\kappa_0=-0.0065, \quad  \kappa_1=-5, \quad  \kappa_2=20.
\end{equation}
The parameters of inflation are
\begin{equation}
	\eta_p = 1.4, \quad \xi=1.2, \quad \gamma=0.8.
\end{equation}
It will be observed from the figures that these parameters, along with parameter $r$, have major effect on stability.

Corresponding to these parameters, we find two equilibrium points 
\begin{align}
	&E_{4,1}(0.836260,0.968365,0.063277), \nonumber \\
	&E_{4,2}(0.808446,0.965992,0.990423).
\end{align}
When $\tau=0$, for the equilibrium point $E_{4,1}$ the criteria in \eqref{RH} are satisfied, therefore it is stable. In this case, $Z(\ws)=0.0049$, therefore the corresponding inflation rate is $0.49\%$. For the second equilibrium point $E_{4,2}$ these criteria do not hold, therefore $E_{4,2}$ is unstable when $\tau=0$. For this case $Z(\ws)=-0.0418$, hence the corresponding inflation rate is $-4.18\%$.

Now we shall calculate the critical value $\tau_0$ as suggested by  \eqref{eqtau} and show that $E_{4,1}$ undergoes a Hopf bifurcation at this value.
From \eqref{eqw}, we have the six roots as follows
\begin{equation}
	\mu_{1,...,6} = (\mp 2.157, \mp 1.881, \mp 0.050i).
\end{equation}
In other words, we have four purely imaginary roots $\pm 1.881 i$ and $\pm 2.157i$ for \eqref{deq}.
We find that  $\mu_0=2.157$ and  $\tau_0=0.82998$. The transversality condition is satisfied as the derivative at that point is $h'(\mu_{0}^{2}) = 5.180 > 0$.

With the set of parameters, considering Eq. \eqref{A58} of the Appendix section,  for the critical delay parameter $\tau_0 = 0.82998$  we have $c_1(0) = 436.694 + 3390.52i$, which gives $\bar{\mu}_2 < 0$, indicating a subcritical Hopf bifurcation. Since $\beta_2 > 0$, the bifurcating periodic solutions are unstable. $T_2 < 0$ suggests that the period of the bifurcating periodic solutions decreases.

We observe by numerical trial that  for the values $\tau<\tau_0$ close to $\tau_0$, the real part of all the eigenvalues found are  negative and for $\tau>\tau_0$ close to  $\tau_0$ there is one eigenvalue for which the real part is positive, exhibiting  the transition from stability to instability.
In Figure 1, for $\tau=0$ the system is shown as stable around the equilibrium point $E_{4,1}$. Figure 2 is at the critical value $\tau=\tau_0=0.82998$. In Figure 3 the plot is  for $\tau=0.85>\tau_0$, and for this value of $\tau$, there is an eigenvalue $0.00579485+2.15435i$ with positive real part. These demonstrate numerically that the system undergoes a Hopf bifurcation at $\tau=\tau_0=0.82998$.

\section{Conclusion}

In this work we have proposed and analysed  a new macroeconomic dynamics. Based on the system constructed by \cite{Grasselli2015}, we have established  the new system \eqref{main}, in which the inflation function makes a delayed entry to the dynamics. We analyzed this system's dynamical behaviors, stability and Hopf bifurcation at the equilibrium points. For any value of the inflation rate at a certain time,  the state variables wage share, employment rate and debt ratio are affected from this value after the delay time $\tau$. Our analysis shows that, an equilibrium point which is asymptotically stable in the delay-free context may become unstable and experience periodic oscillations, therefore a Hopf bifurcation occurs if the delay time exceeds a certain value.

To the best of our knowledge, this is the first time that a bifurcation analysis of the type we have done in this work was applied on a macroeconomic model including the economic variables we have considered here. In the available literature this work was motivated, there are extensions of the Keen and Goodwin's models with further considerations. We believe this work will be an initial spark for further research in this direction.

\section*{Declarations}

{\bf Funding:} This work was supported by Scientific Research Projects Department of Istanbul Technical University, Project Number: PTA-2023-44984. 

\noindent {\bf Credit Author Statement:} All authors contributed to the study equally. S. Harman has the additional role of supervision, C. Özemir has the additional roles of supervision, funding acquisition and project administration.

\appendix
\section{Appendix: Direction of Hopf bifurcation}
By applying the transformations $t=\tau \tilde t $, $\tilde \omega(\tilde t )=\omega(\tau \tilde t)$, $\tilde \lambda(\tilde t )=\lambda(\tau \tilde t)$, 
$\tilde b(\tilde t )=b(\tau \tilde t)$ and dropping the tildes afterwards, we obtain from \eqref{main} the system
\begin{equation}
	\label{apmain}
	\begin{aligned}
		\dot{\om}(t) &= \tau \Bigg[\frac{\Phi_1}{(1-\lambda)^2}-\Phi_0-\alpha-(1-\gamma)\Big(\eta_p (\xi \omega(t-1)-1)\Big)\Bigg]\om(t), \\
		\dot{\lambda}(t) &= \tau \big[g(\pi) - \alpha - \beta\big] \lambda(t), \\
		\dot{b}(t) &= \tau \Big\{\kappa(\pi) - \pi -  \big[\eta_p \xi \om(t-1) + g(\pi) - \eta_p \big] b(t)\Big\}.
	\end{aligned}
\end{equation}
In Section 2.2 we showed that the system \eqref{main}  undergoes a Hopf bifurcation at $E=E(\om^*, \lambda^*, b^*).$
Let $u_1(t)=\om(t)-\om^*$,   $u_2(t)=\lambda(t)-\lambda^*$,  $u_3(t)=b(t)-b^*$, which translates the equilibrium point to the origin, and let $\tau=\tau_k+\chi$. So we obtain 
\begin{equation}
	\label{aptmain}
	\begin{aligned}
		{\dot u_1}(t) &= \tau \Bigg[\frac{\Phi_1}{(1-\lambda^*-u_2(t))^2}+\hat a_0 u_1(t-1)+ \hat a_0 \om^*+\hat a_1\Bigg] [u_1(t) +\om^*], \\
		{\dot  u_2}(t) &= \tau \big[g(\pi_s) - \alpha - \beta \big] [u_2(t)+\lambda^*], \\
		{\dot u_3}(t) &= \tau \Big\{\kappa(\pi_s) - \pi_s -  \big[\eta_p \xi u_1(t-1) + g(\pi_s) +\eta_p \xi \om^*- \eta_p \big] [u_3(t)+b^*]\Big\}
	\end{aligned}
\end{equation}
where $\pi_s=\pi^*-u_1(t)-r u_3(t)$, $\hat a_0 = \xi \eta_p (\gamma-1)$ and 
$\hat a_1 = -\Phi_0-\alpha-\eta_p-\gamma \eta_p $.
After employing the series expansions, we get

\begin{align}
		{\dot u_1}(t) &= \tau\Big[ \frac{2\Phi_1 \om^*}{(1-\lambda^*)^3} u_2(t) +\om^* \hat a_0 u_1(t-1)+\frac{1}{2}\frac{6\Phi_1 \om^*}{(1-\lambda^*)^4}u_2^2(t) +\frac{2 \Phi_1}{(1-\lambda^*)^3} u_1(t) u_2(t)   \nonumber\\
                         & +\hat a_0 u_1(t) u_1(t-1)\Big]+\text{h.o.t.}, \nonumber\\
		{\dot u_2}(t) &= \tau\Big[- \lambda^* g'(\pi^*) u_1(t) -  r \lambda^* g'(\pi^*) u_3(t) + \frac{1}{2} \lambda^* g''(\pi^*)  u_1^2(t) + \frac{1}{2} \lambda^* r^2 g''(\pi^*)  u_3^2(t)  \nonumber\\&
  -  g'(\pi^*) u_1(t) u_2(t)+  r \lambda^* g''(\pi^*) u_1(t) u_3(t)-  r g'(\pi^*) u_2(t) u_3(t)\Big]+\text{h.o.t.}, \nonumber\\
		{\dot u_3}(t) &= \tau \Bigg[1- \kappa'(\pi^*)+b^*g'(\pi^*)\Bigg] u_1(t) + \tau \Bigg[r[1-\kappa'(\pi^*)+b^*g'(\pi^*) ]-[Z(\om^*)+g(\pi^*)]\Bigg] u_3(t)  \label{linearmain}\\&
 - \tau b^* \eta_p \xi u_1(t-1) + \frac{\tau}{2} \Bigg[\kappa''(\pi^*)-g''(\pi^*) b^* \Bigg] u_1^2(t) \nonumber \\
 &+ \frac{\tau}{2} \Bigg[r^2 \kappa''(\pi^*)+2 r g'(\pi^*)-r^2 g''(\pi^*) b^* \Bigg] u_3^2(t) \nonumber \\ &
  + \tau \Bigg[r \kappa''(\pi^*) - r g''(\pi^*) b^* + g'(\pi^*) \Bigg] u_1(t) u_3(t) - \tau \eta_p \xi u_3(t) u_1(t-1)+ \text{h.o.t.}    \nonumber
\end{align}
For the sake of clarity, in the next part of the paper, we will use the following expressions, with $K_0$ to $K_7$ having been introduced in the previous section as

\begin{align}
  K_0&=(\gamma-1) \eta_p \xi \om^*,  &
  K_1&=\om^* \Phi'(\lambda^*), \nonumber\\
  K_2&=\frac{\lambda^* \kappa'(\pi^*)}{\nu},  &
  K_3&=1- \kappa'(\pi^*)+b^* g'(\pi^*),  \nonumber\\
  K_4&= rK_3-[Z(\om^*)+\alpha+\beta],   &
  K_5&=Z(\om^*)+\alpha+\beta-b^* r \eta_p \xi,  \nonumber\\
  K_6&=-b^* \eta_p \xi,     &
  K_7&=Z(\om^*)+\alpha+\beta,\\
  K_8&= \lambda^* \frac{\kappa''(\pi^*)}{2 \nu},  &
  K_9&=  \frac{\kappa''(\pi^*)}{2 \nu} (\nu-b^*),  \nonumber\\
  K_{10}&= r^2 K_9+r \frac{\kappa'(\pi^*)}{\nu},  &
  K_{11}&=2 r K_9+\frac{\kappa'(\pi^*)}{\nu}.  \nonumber
\end{align}
We can write \eqref{linearmain} in the matrix form as 
\begin{equation}
\begin{aligned}
 \begin{bmatrix}
     {\dot u_1}(t)\\
     { \dot u_2}(t)\\
     {\dot u_3}(t)
 \end{bmatrix}&=\tau \begin{bmatrix}
     0 & K_1 & 0\\
     -K_2 & 0 & -r K_2\\
     K_3 & 0 & K_4
 \end{bmatrix} \begin{bmatrix}
     u_1(t)\\
     u_2(t)\\
     u_3(t)
 \end{bmatrix} + \tau \begin{bmatrix}
     K_0 & 0 & 0\\
     0 & 0 & 0\\
     K_6 & 0 & 0
 \end{bmatrix} \begin{bmatrix}
     u_1(t-1)\\
     u_2(t-1)\\
     u_3(t-1)
 \end{bmatrix} \\
 &+ \tau \begin{bmatrix}
      \frac{1}{2} \Phi''(\lambda^*) u_2^2(t)+\Phi'(\lambda^*) u_1(t) u_2(t)+\hat{a}_0 u_1(t) u_1(t-1)\\
         K_8 u_1^2(t)+r^2 K_8 u_3^2(t)-g'(\pi^*) u_1(t) u_2(t)+r \lambda^* g''(\pi^*) u_1(t) u_3(t)-r g'(\pi^*) u_2(t) u_3(t)\\
         K_9 u_1^2(t)+K_{10} u_3^2(t)+K_{11} u_1(t) u_3(t)-\eta_p \xi u_3(t) u_1(t-1)
 \end{bmatrix}
 \end{aligned}
\end{equation}
including terms up to second order.
Let us define
\begin{equation}
    u_t(\theta)=\begin{bmatrix}
        u_{1t}(\theta)\\
        u_{2t}(\theta)\\
        u_{3t}(\theta)
    \end{bmatrix}:=u(t+\theta)=\begin{bmatrix}
        u_1(t+\theta)\\
         u_2(t+\theta)\\
          u_3(t+\theta)
    \end{bmatrix}, \quad u_t:[-1,0] \rightarrow \mathbb{R}^3.
\end{equation}
Therefore we get
\begin{equation}
\begin{aligned}
 &\begin{bmatrix}
     {\dot u_1}(t)\\
     {\dot u_2}(t)\\
      {\dot u_3}(t)
 \end{bmatrix}=\tau \begin{bmatrix}
     0 & K_1 & 0\\
     -K_2 & 0 & -r K_2\\
     K_3 & 0 & K_4
 \end{bmatrix} \begin{bmatrix}
     u_{1t}(0)\\
     u_{2t}(0)\\
     u_{3t}(0)
 \end{bmatrix} + \tau \begin{bmatrix}
      K_0 & 0 & 0\\
     0 & 0 & 0\\
     K_6 & 0 & 0
 \end{bmatrix} \begin{bmatrix}
     u_{1t}(-1)\\
     u_{2t}(-1)\\
     u_{3t}(-1)
 \end{bmatrix} \\
  &+\tau \begin{bmatrix}
      \frac{1}{2} \Phi''(\lambda^*) u_{2t}^2(0)+\Phi'(\lambda^*) u_{1t}(0) u_{2t}(0)+\hat{a}_0 u_{1t}(0) u_{1t}(-1)\\
         K_8 u_{1t}^2(0)+r^2 K_8 u_{3t}^2(0)-g'(\pi^*) u_{1t}(0) u_{2t}(0)+r \lambda^* g''(\pi^*) u_{1t}(0) u_{3t}(0)-r g'(\pi^*) u_{2t}(0) u_{3t}(0)\\
    K_9 u_{1t}^2(0)+K_{10} u_{3t}^2(0)+K_{11} u_{1t}(0) u_{3t}(0)-\eta_p \xi u_{3t}(0) u_{1t}(-1)
 \end{bmatrix}
 \end{aligned}
\end{equation}
For $\phi(\theta)=(\phi_1(\theta),\phi_2(\theta), \phi_3(\theta)) \in C([-1,0],\mathbb{R}^3)$, define the operator $L_\chi$ by

\begin{equation}
    L_\chi=\tau \begin{bmatrix}
     0 & K_1 & 0\\
     -K_2 & 0 & -r K_2\\
     K_3 & 0 & K_4
    \end{bmatrix} \begin{bmatrix}
        \phi_1(0)\\
        \phi_2(0)\\
        \phi_3(0)
    \end{bmatrix}+ \tau \begin{bmatrix}
         K_0 & 0 & 0\\
     0 & 0 & 0\\
     K_6 & 0 & 0
    \end{bmatrix} \begin{bmatrix}
         \phi_1(-1)\\
        \phi_2(-1)\\
        \phi_3(-1)
    \end{bmatrix}
\end{equation}
and
\begin{equation}
\begin{aligned}
    f(\phi, \chi)
    &= \tau \begin{bmatrix}
         \frac{1}{2} \Phi''(\lambda^*) \phi_2^2(0)+\Phi'(\lambda^*) \phi_1(0) \phi_2(0)+\hat{a}_0 \phi_1(0) \phi_1(-1)\\
         K_8 \phi_1^2(0)+r^2 K_8 \phi_3^2(0)-g'(\pi^*) \phi_1(0) \phi_2(0)+r \lambda^* g''(\pi^*) \phi_1(0) \phi_3(0)-r g'(\pi^*) \phi_2(0) \phi_3(0)\\
    K_9 \phi_1^2(0)+K_{10} \phi_3^2(0)+K_{11} \phi_1(0) \phi_3(0)-\eta_p \xi \phi_3(0) \phi_1(-1)
    \end{bmatrix}\\
    &:=\tau \begin{bmatrix}
        f_{11}\\
        f_{12}\\
        f_{13} 
    \end{bmatrix}
    \end{aligned}
\end{equation}
Therefore we have
\begin{equation}
    \dot{u}(t)=L_{\chi}(u_t)+f(u_t,\chi).
\end{equation}
By the Riesz representation theorem, there exists a matrix whose components are bounded variation functions $\eta(\theta,\chi)$ in $\theta \in [-1,0]$ such that
\begin{equation}
     L_{\chi} \phi = \int_{-1}^{0} d\eta(\theta,\chi) \phi(\theta)
\end{equation}
for $\phi \in C[-1,0]$. $\eta(\theta,\chi)$ can be chosen as
\begin{equation}
    \eta(\theta,\chi)=\tau \begin{bmatrix}
        0 & K_1 & 0\\
        -K_2 & 0 & -r K_2\\
        K_3 & 0 & K_4
    \end{bmatrix} \delta(\theta)+ \tau \begin{bmatrix}
        K_0 & 0 & 0\\
        0 & 0 & 0\\
        K_6 & 0 & 0
    \end{bmatrix} \delta(\theta+1)
\end{equation}
where $\delta(\theta)$ is the Dirac distribution.
For $\phi(\theta) \in C^1([-1,0], \mathbb{R}^3)$, define
\begin{equation}
    A(\chi) \phi(\theta)= \begin{cases}
           \quad \dfrac{d\phi(\theta)}{d\theta}, & \theta \in [-1,0),\\
          \displaystyle \int_{-1}^{0} d\eta(s,\chi) \phi(s), & \theta=0,
      \end{cases}
\end{equation}
and
\begin{equation}
    R(\chi) \phi(\theta)= \begin{cases}
           0, & \theta \in [-1,0),\\
           f(\phi,\chi), & \theta=0,
      \end{cases}
\end{equation}
So we obtain that the system under consideration is equivalent to 
\begin{equation} \label{eqmain}
    \dot{u}(t) = A(\chi) u_t + R(\chi) u_t.
\end{equation}
For $\psi \in C^1([0,1],(\mathbb{R}^3)^*)$, define the adjoint $A^*(0)$ of $A(0)$ as
\begin{equation}
    A^*(0) \psi(s)= \begin{cases}
           -\dfrac{d\psi(s)}{ds}, & s \in (0,1],\\
           \displaystyle\int_{-1}^{0} d\eta^T(t,0) \psi(-t), & s=0,
      \end{cases}
\end{equation}
where $T$ denotes transpose and a bilinear product
\begin{equation}
\label{bin}
    \langle \psi, \phi \rangle = \Bar{\psi}(0) \cdot \phi(0)- \int_{\theta=-1}^{0} \int_{\xi=0}^{\theta} \Bar{\psi}^T (\xi-\theta) d\eta(\theta) \phi(\xi) d\xi,
\end{equation}
here $\eta(\theta)=\eta(\theta,0)$ and a bar denotes complex conjugate. Then $A(0)$ and $A^*(0)$ are adjoint operators. Also we know that $\pm i \tau_k \mu$ are eigenvalues of $A(0)$. Thus, they are also eigenvalues of $A^*$. Let $q(\theta)$ be the eigenvector of $A(0)$ corresponding to $i \tau_k \mu$ and $q^*(s)$ be the eigenvector of $A^*(0)$ corresponding to $- i \tau_k \mu.$ In this case, we obtain these eigenvectors as 
\begin{equation}
    q(\theta)=\begin{bmatrix}
        1\\
        \dfrac{-K_0 e^{-i \mu \tau_k}+i \mu}{K_1}\\
        \dfrac{-K_3-K_6 e^{-i \mu \tau_k}}{K_4 - i\mu}
    \end{bmatrix} e^{i \mu \tau_k \theta}:=\begin{bmatrix}
        1\\
        \alpha\\
        \beta
    \end{bmatrix} e^{i \mu \tau_k \theta}
\end{equation}
and
\begin{equation}
    q^*(s)=B \begin{bmatrix}
        1\\
       \dfrac{- K_1}{i \mu}\\
        \dfrac{- r K_1 K_2}{i \mu K_4 - \mu^2}
    \end{bmatrix} e^{i \mu \tau_k s}:= B \begin{bmatrix}
        1\\
        \alpha^*\\ 
        \beta^*
    \end{bmatrix} e^{i \mu \tau_k s}
\end{equation}
By applying \eqref{bin} inner product formula, we get
\begin{equation}
\begin{aligned}
    \langle q^*(s),q(\theta) \rangle &=\bar q^*(0)^T q(0)-\int_{\theta=-1}^{0} \int_{\xi=0}^{\theta} \bar q^*(\xi-\theta)^T d\eta(\theta) q(\xi) d\xi\\
    &=\bar B \begin{bmatrix}
         1 & \bar \alpha^* &  \bar \beta^*
    \end{bmatrix} \begin{bmatrix}
        1\\
        \alpha \\
        \beta
    \end{bmatrix}
    -\int_{\theta=-1}^{0} \int_{\xi=0}^{\theta} \bar B \begin{bmatrix}
          1 & \bar \alpha^* & \bar \beta^*
    \end{bmatrix} e^{-i \mu \tau_k (\xi-\theta)} d\eta(\theta) \begin{bmatrix}
        1\\
        \alpha \\
        \beta
    \end{bmatrix} e^{i \mu \tau_k \xi} d\xi\\
    &= \bar B \Biggl\{1+ \bar \alpha^* \alpha + \bar \beta^* \beta-\int_{\theta=-1}^{0} \begin{bmatrix}
         1 & \bar \alpha^* & \bar \beta^*
    \end{bmatrix} e^{i \mu \tau_k \theta} \theta d \eta(\theta) \begin{bmatrix}
         1\\
        \alpha \\
        \beta
    \end{bmatrix} \Biggr\}\\
    &=\bar B ( 1+ \bar \alpha^* \alpha + \bar \beta^* \beta ) \\
    &-\int_{\theta=-1}^{0} \bar B \begin{bmatrix}
         1 & \bar \alpha^* & \bar \beta^*
    \end{bmatrix} e^{i \mu \tau_k \theta}\theta 
     \Biggl\{\tau_k \begin{bmatrix}
       0 & K_1 & 0\\
        -K_2 & 0 & -r K_2\\
        K_3 & 0 & K_4
    \end{bmatrix} \delta(\theta) \\
    & \qquad \qquad \qquad \qquad \qquad \qquad \qquad +  \tau_k \begin{bmatrix}
        K_0 & 0 & 0\\
        0 & 0 & 0\\
        K_6 & 0 & 0
    \end{bmatrix} \delta(\theta+1) \Biggr\} \begin{bmatrix}
         1\\
        \alpha \\
        \beta
    \end{bmatrix} d\theta\\
    &=\bar B ( 1+ \bar \alpha^* \alpha + \bar \beta^* \beta ) + \bar B \begin{bmatrix}
        1 & \bar \alpha^* & \bar \beta^*
    \end{bmatrix}\begin{bmatrix}
        \tau_k K_0 & 0 & 0\\
        0 & 0 & 0\\
        \tau_k K_6 & 0 & 0
    \end{bmatrix}  \begin{bmatrix}
         1\\
        \alpha \\
        \beta
    \end{bmatrix} e^{-i \mu \tau_k}\\
    &=\bar B \Big( 1+\bar \alpha^* \alpha + \bar \beta^* \beta + \tau_k ( K_0 + \bar \beta^* K_6 )e^{-i \mu \tau_k} \Big)
\end{aligned}
\end{equation}
Then we choose
\begin{equation}
    \bar B= \frac{1}{1+\bar \alpha^* \alpha + \bar \beta^* \beta + \tau_k ( K_0 + \bar \beta^* K_6 )e^{-i \mu \tau_k}}
\end{equation}
which assures that 
\begin{equation}
    \langle q^*(s), q(\theta) \rangle =1.
\end{equation}

Now we need to calculate the coordinates to define the center manifold $\bold C_0$ at $\chi=0.$\\
Let $u_t$ be the solution of Eq. \eqref{linearmain} when $\chi=0.$ Define
\begin{equation}\label{defzw}
    z(t)=\langle q^*,u_t \rangle , \qquad  W(t,\theta)=u_t(\theta)-2  \text{Re} \{z(t) q(\theta)\}.
\end{equation}
On the center manifold $\bold C_0$
\begin{equation}\label{defw}
    W(t,\theta)=W(z(t),\bar z(t),\theta),
\end{equation}
where
\begin{equation}\label{new315}
    W(z,\bar z,\theta)=W_{20}(\theta) \frac{z^2}{2}+W_{11}(\theta) z \bar z+ W_{02}(\theta) \frac{\bar z^2}{2}+W_{30} (\theta) \frac{z^3}{6}+...
\end{equation}
Here $z$ and $\bar z$ are local coordinates for center manifold $\bold C_0$ in the direction of $q^*$ and $\bar q^*.$ It can be easily seen that if $u_t$ is real then $W$ is real. In this study only real solutions are considered. For the solution $u_t \in \bold C_0$ of \eqref{linearmain}, since $\chi=0$, using inner product properties we obtain
\begin{equation}\label{new316}
    \dot z(t)=i \mu\tau_k z+ \bar q^*(\theta)^T f(0,u_t)=i \mu \tau_k z+ \bar q^*(0)^T f_0(z, \bar z)=i \mu \tau_k z+g(z,\bar z),
\end{equation}
where $\bar q^*(0)^T f_0(z, \bar z)=g(z,\bar z)$ and 
\begin{equation}\label{gz}
    g(z,\bar z)=g_{20} \frac{z^2}{2}+g_{11} z \bar z+ g_{02} \frac{\bar z^2}{2}+g_{21} \frac{z^2 \bar z}{2}+\dots 
\end{equation}
Let us note that

\begin{equation}
    u_t(u_{1t}(\theta), u_{2t}(\theta), u_{3t}(\theta))=W(t,\theta)+z q(\theta)+ \bar z \bar q(\theta),
\end{equation}
and $ q(\theta)= \begin{bmatrix}
    1 & \alpha & \beta
\end{bmatrix} ^T e^{i \mu \tau_k \theta}$, and explicitly,
\begin{equation}
    \begin{aligned}
        u_{1t}(0)&= z+ \bar z + W_{20}^{(1)}(0) \frac{z^2}{2} +  W_{11}^{(1)}(0) z \bar z + W_{02}^{(1)}(0) \frac{\bar z^2}{2}+ ... \\
        u_{2t}(0)&= z \alpha + \bar z \bar \alpha + W_{20}^{(2)}(0) \frac{z^2}{2} +  W_{11}^{(2)}(0) z \bar z + W_{02}^{(2)}(0) \frac{\bar z^2}{2}+ ... \\
        u_{3t}(0)&= z \beta + \bar z \bar \beta + W_{20}^{(3)}(0) \frac{z^2}{2} +  W_{11}^{(3)}(0) z \bar z + W_{02}^{(3)}(0) \frac{\bar z^2}{2}+ ... \\
        u_{1t}(-1)&= z e^{-i \mu \tau_k} + \bar z e^{i \mu \tau_k}+ W_{20}^{(1)}(-1) \frac{z^2}{2} + W_{11}^{(1)}(-1) z \bar z + W_{02}^{(1)}(-1) \frac{\bar z^2}{2}+ ... \\
        u_{2t}(-1)&= z \alpha e^{-i \mu \tau_k} + \bar z \bar \alpha e^{i \mu \tau_k}+ W_{20}^{(2)}(-1) \frac{z^2}{2} +  W_{11}^{(2)}(-1) z \bar z + W_{02}^{(2)}(-1) \frac{\bar z^2}{2}+ ...\\
        u_{3t}(-1)&= z \beta e^{-i \mu \tau_k} + \bar z \bar \beta e^{i \mu \tau_k}+ W_{20}^{(3)}(-1) \frac{z^2}{2} + W_{11}^{(3)}(-1) z \bar z + W_{02}^{(3)}(-1) \frac{\bar z^2}{2}+ ...\\
    \end{aligned}
\end{equation}

and 
\begin{equation} \label{gzf}
     g(z, \bar z) = \bar q^*(0)^T f_0(z, \bar z) = \bar B \begin{bmatrix}
         1 & \bar \alpha^* & \bar \beta^*
     \end{bmatrix} \tau_k \begin{bmatrix}
         f_{11}^0 \\
         f_{12}^0\\
         f_{13}^0
     \end{bmatrix}\\
\end{equation}
where 
\begin{equation}
    \begin{bmatrix}
         f_{11}^0 \\
         f_{12}^0\\
         f_{13}^0
     \end{bmatrix}=\begin{bmatrix}
      \frac{1}{2} \Phi''(\lambda^*) u_{2t}^2(0)+\Phi'(\lambda^*) u_{1t}(0) u_{2t}(0)+\hat{a}_0 u_{1t}(0) u_{1t}(-1)\\
         K_8 u_{1t}^2(0)+r^2 K_8 u_{3t}^2(0)-g'(\pi^*) u_{1t}(0) u_{2t}(0)+r \lambda^* g''(\pi^*) u_{1t}(0) u_{3t}(0)-r g'(\pi^*) u_{2t}(0) u_{3t}(0)\\
    K_9 u_{1t}^2(0)+K_{10} u_{3t}^2(0)+K_{11} u_{1t}(0) u_{3t}(0)-\eta_p \xi u_{3t}(0) u_{1t}(-1)
 \end{bmatrix}.
\end{equation}
Substituting \eqref{gz} into the left side of \eqref{gzf} and comparing the coefficients, we obtain 
\begin{equation}
\begin{aligned}
    g_{20}&=2 \bar B \tau_k \Bigl ( \frac{\alpha^2}{2} \Phi''(\lambda^*)+ \alpha \Phi'(\lambda)+\hat{a}_0 e^{- i \mu \tau_k}+\bar \alpha^* K_8 + \bar \alpha^* \beta^2 r^2 K_8 -\alpha \bar \alpha^* g'(\pi^*)\\
    &+ \beta \bar \alpha^* r \lambda^* g''(\pi^*) - \alpha \bar \alpha^* \beta r g'(\pi^*) + \bar \beta^* K_9+ \beta^2 \bar \beta^* K_{10} + \beta \bar \beta^* K_{11} - \beta \bar \beta^* \eta_p \xi e^{-i \mu \tau_k}\Bigr),\\
    g_{11}&= 2 \bar B \tau_k \Bigl( \frac{\alpha \bar \alpha}{2} \Phi''(\lambda^*)+ \text{Re}\{\alpha\}\Phi'(\lambda^*) + \text{Re}\{e^{i \mu \tau_k}\} \hat{a}_0+\bar \alpha^* K_8+ \beta \bar \beta \bar \alpha^* r^2 K_8 \\
    &- \bar \alpha^* \text{Re}\{ \alpha\} g'(\pi^*)+ \bar \alpha^* \text{Re}\{\beta\} r \lambda^* g''(\pi^*) - \bar \alpha^* \text{Re}\{ \alpha \beta \} r g'{\pi^*}+ \bar \beta^* K_9+ \beta \bar \beta \bar \beta^* K_{10} \\
    &  + \bar \beta^* \text{Re} \{ \beta \} K_{11} - 2 \bar \beta^* \eta_p \xi \text{Re}\{\beta e^{i \mu \tau_k} \} \Bigr), \\
    g_{02}&= 2 \bar B \tau_k \Bigl(  \frac{\bar \alpha^2}{2} \Phi''(\lambda^*)+ \bar \alpha \Phi'(\lambda)+\hat{a}_0 e^{i \mu \tau_k}+\bar \alpha^* K_8 + \bar \alpha^* \bar \beta^2 r^2 K_8 -\bar \alpha \bar \alpha^* g'(\pi^*)\\
    &+ \bar \beta \bar \alpha^* r \lambda^* g''(\pi^*) - \bar \alpha \bar \alpha^* \bar \beta r g'(\pi^*) + \bar \beta^* K_9+\bar \beta^2 \bar \beta^* K_{10} +\bar \beta \bar \beta^* K_{11} - \bar \beta \bar \beta^* \eta_p \xi e^{i \mu \tau_k} \Bigr),            \\
    g_{21}&=\bar B \tau_k \Biggl( \Phi''(\lambda^*)\Bigl( \bar \alpha W_{20}^{(2)} (0)+2 \alpha W_{11}^{(2)}(0) \Bigr) \\
    &+ \Phi'(\lambda^*) \Bigl( \bar \alpha W_{20}^{(1)}(0)+W_{20}^{(2)}(0)+2 \alpha W_{11}^{(1)}(0)+2 W_{11}^{(2)}(0) \Bigr)\\
    &+\hat{a}_0 \Bigl(2 W_{11}^{(1)}(-1) +W_{20}^{(1)}(-1)+ e^{i \mu \tau_k} W_{20}^{(1)}(0)+ 2 e^{- i \mu \tau_k} W_{11}^{(1)}(0) \Bigr)\\ 
    &+ \bar \alpha^* K_8 \Bigl(2 W_{20}^{(1)}(0) +4  W_{11}^{(1)}(0)\Bigr)
    +\bar \alpha^* r^2 K_8 \Bigl( 2 \bar \beta W_{20}^{(3)}(0)+4 \beta W_{11}^{(3)}(0) \Bigr)\\
    &- \bar \alpha^* g'(\pi^*) \Bigl( \bar \alpha W_{20}^{(1)}(0)+W_{20}^{(1)}(0)+ 2 \alpha W_{11}^{(1)}(0)+2 W_{20}^{(1)}(0) \Bigr) \\
    &+\bar \alpha^* \lambda^* r g''(\pi^*) \Bigl(W_{20}^{(1)}(0)+W_{20}^{(1)}(0)+W_{20}^{(1)}(0)+W_{11}^{(2)}(0) \Bigr)\\
    &- \bar \alpha^* r g'(\pi^*) \Bigl( \bar \beta W_{20}^{(2)}(0)+ \alpha W_{20}^{(3)}(0)+2 \beta W_{11}^{(3)}(0)+2 \alpha W_{11}^{(3)}(0) \Bigr) \\
    &+\bar \beta^* K_9 \Bigl( 2 W_{20}^{(1)}(0)+4 W_{11}^{(1)}(0)\Bigr) +\bar \beta^* K_{10} \Bigl(2 \bar \beta W_{20}^{(3)}(0)+4 \beta W_{11}^{(3)}(0) \Bigr) \\
    &+ \bar \beta^* K_{11}\Bigl(\bar \beta W_{20}^{(1)}(0)+W_{20}^{(3)}(0)+2 \beta W_{11}^{(1)}(0)+2 W_{11}^{(3)}(0) \Bigr) \\
    &- \bar \beta^* \eta_p \xi \Bigl( \bar \beta W_{20}^{(1)}(-1)+e^{i \mu \tau_k }W_{20}^{(3)}(0)+2 \beta W_{11}^{(1)}(-1)+2 e^{-i \mu \tau_k} W_{11}^{(3)}(0) \Bigr) \Biggr).
\end{aligned}
\end{equation}
    We still need to compute $W_{20}(\theta)$ and  $W_{11}(\theta)$ to find the value of the coefficient $g_{21}.$
  From \eqref{eqmain} and \eqref{defzw}, we have
\begin{align}\label{new322}
     \begin{split}
        \Dot W= \dot u_t-\dot z q - \dot{\bar z }\bar q=& \begin{cases}
            A W - 2 \text{Re}\{ \bar q^*(0) f_0 q(\theta)\}, \qquad \theta \in [-1,0), \\
             A W - 2 \text{Re}\{ \bar q^*(0) f_0 q(0)\}+ f_0, \qquad \theta=0,
        \end{cases}\\
        =& A W + H(z, \bar z, \theta),
    \end{split}
\end{align}
   where
   \begin{equation}\label{Hz}
       H(z, \bar z, \theta)=H_{20}(\theta) \frac{z^2}{2} + H_{11}(\theta) z \bar z + H_{02}(\theta) \frac{\bar z^2}{2}+ \dots
   \end{equation}
If we substitute the corresponding series in \eqref{new322} and compare the coefficients, we obtain
\begin{equation}\label{match}
    (A(0)- 2 i \mu \tau_k ) W_{20}(\theta) = -H_{20}(\theta), \qquad  A(0) W_{11}(\theta)=-H_{11}(\theta).
\end{equation}
Also from \eqref{gz}, we can write for $ \theta \in [-1,0)$,
\begin{equation}\label{H2}
    H(z, \bar z, \theta)=- \bar q^*(0)^Tf_0 q(\theta) - q^*(0)^T \bar f_0 \bar q(\theta)=-g(z,\bar z) q(\theta) - \bar g(z,\bar z) \bar q(\theta).
\end{equation}
Then, by comparing the coefficients with \eqref{Hz}, we determine that
\begin{equation}\label{H20}
    H_{20}(\theta) = -g_{20} q(\theta) - \bar g_{02} \bar q(\theta)
\end{equation}
and
\begin{equation}\label{H11}
     H_{11}(\theta) = -g_{11} q(\theta) - \bar g_{11} \bar q(\theta).
\end{equation}
From \eqref{match} and \eqref{H20}, we obtain
\begin{equation*}
    \dot W_{20}(\theta)=2 i \mu \tau_k W_{20}(\theta)+ g_{20} q(\theta)+ \bar g_{02} \bar q(\theta).
\end{equation*}
If we solve the last equation taking into account $q(\theta)=q(0) e^{i \mu \tau_k \theta}$, we obtain the following
\begin{equation}\label{W20}
    W_{20}(\theta)=\frac{i g_{20}}{\mu \tau_k} q(0) e^{i \mu \tau_k \theta}+ \frac{i \bar g_{02}}{3 \mu \tau_k} \bar q(0) e^{- i \mu \tau_k \theta}+ E_1 e^{2 i \mu \tau_k \theta}.
\end{equation}
Similarly, from \eqref{match} and \eqref{H11}, we get
\begin{equation}
    \Dot W_{11}(\theta)=g_{11} q(\theta)+ \bar g_{11} \bar q(\theta)
\end{equation}
 and
\begin{equation}\label{W11}
    W_{11}(\theta)=-\frac{i g_{11}}{\mu \tau_k} q(0) e^{i \mu \tau_k \theta}+ \frac{i \bar g_{11}}{\mu \tau_k} \bar q(0) e^{- i \mu \tau_k \theta}+ E_2.
\end{equation}
We must then determine the corresponding $E_1$ and $E_2$ appearing in \eqref{W20} and \eqref{W11}, respectively. To do this, we can first write the following from the definition of $A$ and from \eqref{defw}
\begin{align}\label{intW20}
    \int_{-1}^{0} d\eta(\theta) W_{20}(\theta) &= 2 i \mu \tau_k W_{20}(\theta)-H_{20}(\theta),\\ \label{intW11}
    \int_{-1}^{0} d\eta(\theta) W_{11}(\theta) &= -H_{11}(\theta). 
\end{align}
Also from \eqref{new322}, we have
\begin{equation}\label{newH20}
     H_{20}(0)= -g_{20} q(0)- \bar g_{02} \bar q(0) + \tau_k \begin{bmatrix}
        \alpha^2 \Phi''(\lambda^*) + 2 \alpha \Phi'(\lambda^*) +2 \hat{a}_0 e^{-i \mu \tau_k} \\
        2 K_8 + 2 r^2 \beta^2 K_8 - 2 \alpha g'(\pi^*) + 2 r \lambda^* \beta  g''(\pi^*) - 2 r \alpha \beta g'(\pi^*) \\
        2 K_9 + 2 \beta^2 K_{10} + 2 \beta K_{11}-2 \beta \eta_p \xi e^{-i \mu \tau_k}
    \end{bmatrix}
\end{equation}
   and
   \begin{equation}\label{newH11}
   \begin{aligned}
         H_{11}(0)&=-g_{11} q(0) - \bar g_{11} \bar q(0) \\
         &+\tau_k \begin{bmatrix}
        \alpha \bar \alpha \Phi''(\lambda^*) + 2 \Phi'(\lambda^*) \text{Re}\{\alpha\}+ 2 \hat{a}_0 \text{Re}\{e^{i \mu \tau_k}\} \\
        2 K_8(1 +  \beta \bar \beta r^2) - 2 g'(\pi^*) \text{Re}\{ \alpha \}+2 r \lambda^* g''(\pi^*) \text{Re} \{ \beta \} -2 r \text{Re}\{ \alpha \beta \} g'(\pi^*) \\
        2 K_9 + 2 \beta \bar \beta K_{10} + 2 \text{Re} \{\beta \} K_{11}-2 \eta_p \xi \text{Re}\{ \beta e^{i \mu \tau_k} \}
    \end{bmatrix}.
   \end{aligned}  
   \end{equation}
Substituting \eqref{W20} and \eqref{newH20} into \eqref{intW20}, we get
\begin{equation}
    \Bigg ( 2 i \mu \tau_k I - \int_{-1}^{0} e^{2 i \mu \tau_k \theta } d\eta(\theta) \Bigg) E_1= \tau_k \begin{bmatrix}
       \alpha^2 \Phi''(\lambda^*) + 2 \alpha \Phi'(\lambda^*) +2 \hat{a}_0 e^{-i \mu \tau_k} \\
        2 K_8 + 2 r^2 \beta^2 K_8 - 2 \alpha g'(\pi^*) + 2 r \lambda^* \beta  g''(\pi^*) - 2 r \alpha \beta g'(\pi^*) \\
        2 K_9 + 2 \beta^2 K_{10} + 2 \beta K_{11}-2 \beta \eta_p \xi e^{-i \mu \tau_k}
    \end{bmatrix},
\end{equation}
that is,
\begin{equation}
\begin{bmatrix}
    2 i \mu - K_0 e^{-2i \mu \tau_k} & -K_1  & 0 \\
    K_2 & 2i \mu & r K_2 \\
    -K_3-K_6 e^{-2i \mu \tau_k} & 0 & 2i \mu - K_4
\end{bmatrix} E_1 = \begin{bmatrix}
  \alpha^2 \Phi''(\lambda^*) + 2 \alpha \Phi'(\lambda^*) +2 \hat{a}_0 e^{-i \mu \tau_k} \\
        2 K_8 + 2 r^2 \beta^2 K_8 - 2 \alpha g'(\pi^*) + 2 r \lambda^* \beta  g''(\pi^*) - 2 r \alpha \beta g'(\pi^*) \\
        2 K_9 + 2 \beta^2 K_{10} + 2 \beta K_{11}-2 \beta \eta_p \xi e^{-i \mu \tau_k}
\end{bmatrix}.
\end{equation}
Similarly, substituting \eqref{W11} and \eqref{newH11} into \eqref{intW11}, we obtain
\begin{equation}
    \int_{-1}^{0} d\eta(\theta) E_2= -\tau_k \begin{bmatrix}
        \alpha \bar \alpha \Phi''(\lambda^*) + 2 \Phi'(\lambda^*) \text{Re}\{\alpha\}+ 2 \hat{a}_0 \text{Re}\{e^{i \mu \tau_k}\} \\
        2 K_8(1 +  \beta \bar \beta r^2) - 2 g'(\pi^*) \text{Re}\{ \alpha \}+2 r \lambda^* g''(\pi^*) \text{Re} \{ \beta \} -2 r \text{Re}\{ \alpha \beta \} g'(\pi^*) \\
        2 K_9 + 2 \beta \bar \beta K_{10} + 2 \text{Re} \{\beta \} K_{11}-2 \eta_p \xi \text{Re}\{ \beta e^{i \mu \tau_k} \}
    \end{bmatrix}
\end{equation}
that is
\begin{equation}
\begin{bmatrix}
    K_0 & K_1 & 0 \\
   -K_2 & 0 & -r K_2 \\
   K_3+K_6 & 0 & K_4
\end{bmatrix} E_2 =-\begin{bmatrix}
        \alpha \bar \alpha \Phi''(\lambda^*) + 2 \Phi'(\lambda^*) \text{Re}\{\alpha\}+ 2 \hat{a}_0 \text{Re}\{e^{i \mu \tau_k}\} \\
        2 K_8(1 +  \beta \bar \beta r^2) - 2 g'(\pi^*) \text{Re}\{ \alpha \}+2 r \lambda^* g''(\pi^*) \text{Re} \{ \beta \} -2 r \text{Re}\{ \alpha \beta \} g'(\pi^*) \\
        2 K_9 + 2 \beta \bar \beta K_{10} + 2 \text{Re} \{\beta \} K_{11}-2 \eta_p \xi \text{Re}\{ \beta e^{i \mu \tau_k} \}
    \end{bmatrix}.
\end{equation}
This allows us to write $g_{21}$.  Therefore,  we can now calculate the following values, which characterize the bifurcating periodic solutions at the critical value $\tau_k:$
\begin{align}\label{A58}
    \begin{split}
        c_1(0)=&\frac{i}{2 \mu \tau_k} \Big ( g_{11} g_{20}-2 |g_{11}|^2 - \frac{|g_{02}|^2}{3} \Big ) + \frac{g_{21}}{2},\\
        \bar \mu_2=&-\frac{\text{Re}(c_1(0))}{\text{Re}(x_0'(\tau_k))},\\
        \beta_2=&2 \text{Re}(c_1(0)),\\
        T_2=& - \frac{\text{Im}(c_1(0)) +\bar \mu_2 \text{Im}(x_0'(\tau_k))}{\mu}.
    \end{split}
\end{align}

Here the parameter $\bar{\mu}_2$ determines the direction of the Hopf bifurcation: If $\bar{\mu}_2 > 0$, the Hopf bifurcation is supercritical, and bifurcating periodic solutions exist for $\tau > \tau_k$. If $\bar{\mu}_2 < 0$, the Hopf bifurcation is subcritical, and bifurcating periodic solutions exist for $\tau < \tau_k$.
The stability of the bifurcating periodic solutions is determined by $\beta_2$: If $\beta_2 < 0$, the bifurcating periodic solutions are stable. If $\beta_2 > 0$, the bifurcating periodic solutions are unstable. The period of the bifurcating periodic solutions is governed by $T_2$: If $T_2 > 0$, the period increases. If $T_2 < 0$, the period decreases.

\begin{figure}[h!]
	\centering
	\includegraphics[scale=0.60]{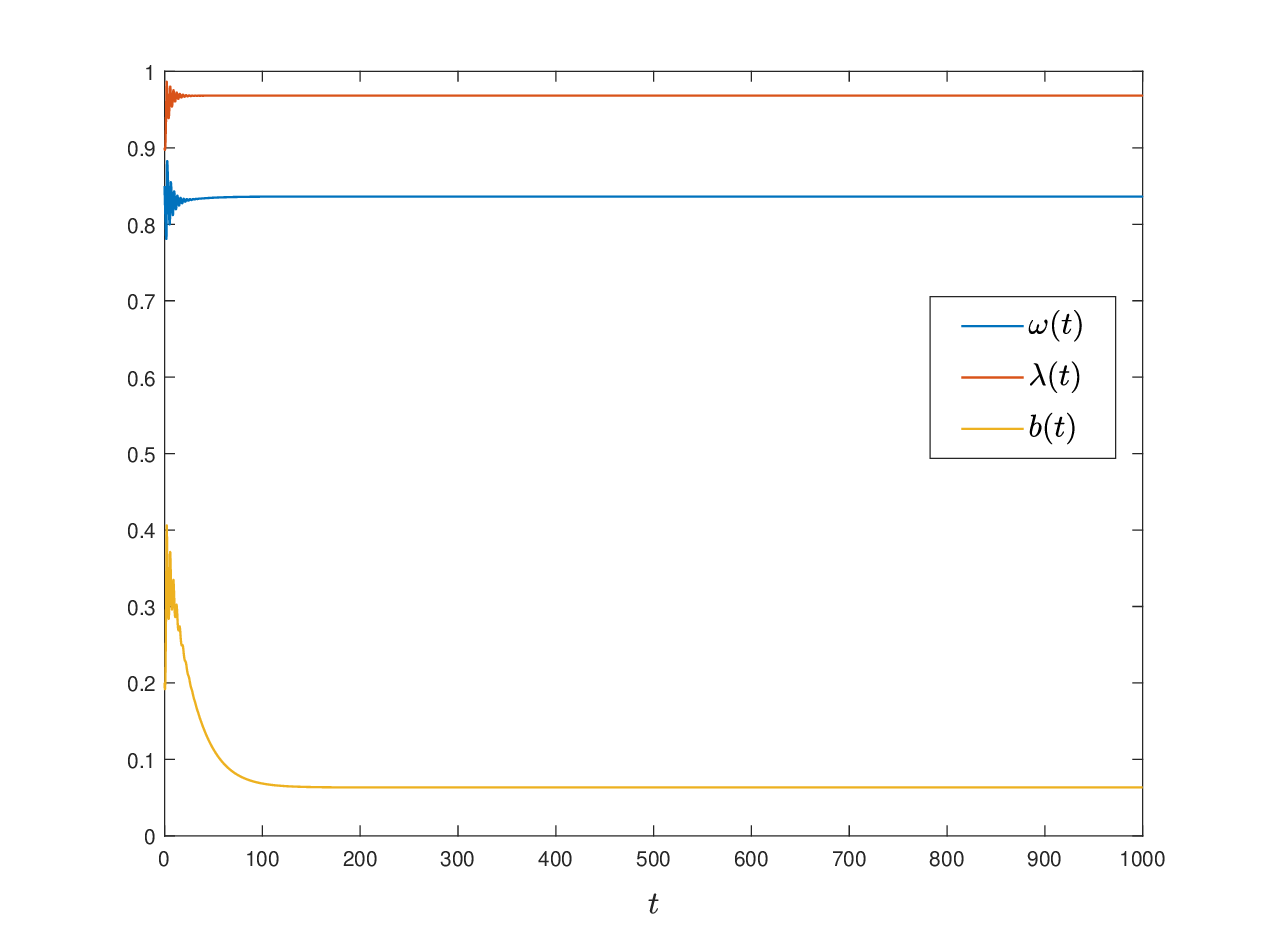}
	\caption{$\tau=0$}
	\label{figure1}
\end{figure}

\begin{figure}[h!]
	\centering
	\includegraphics[scale=0.60]{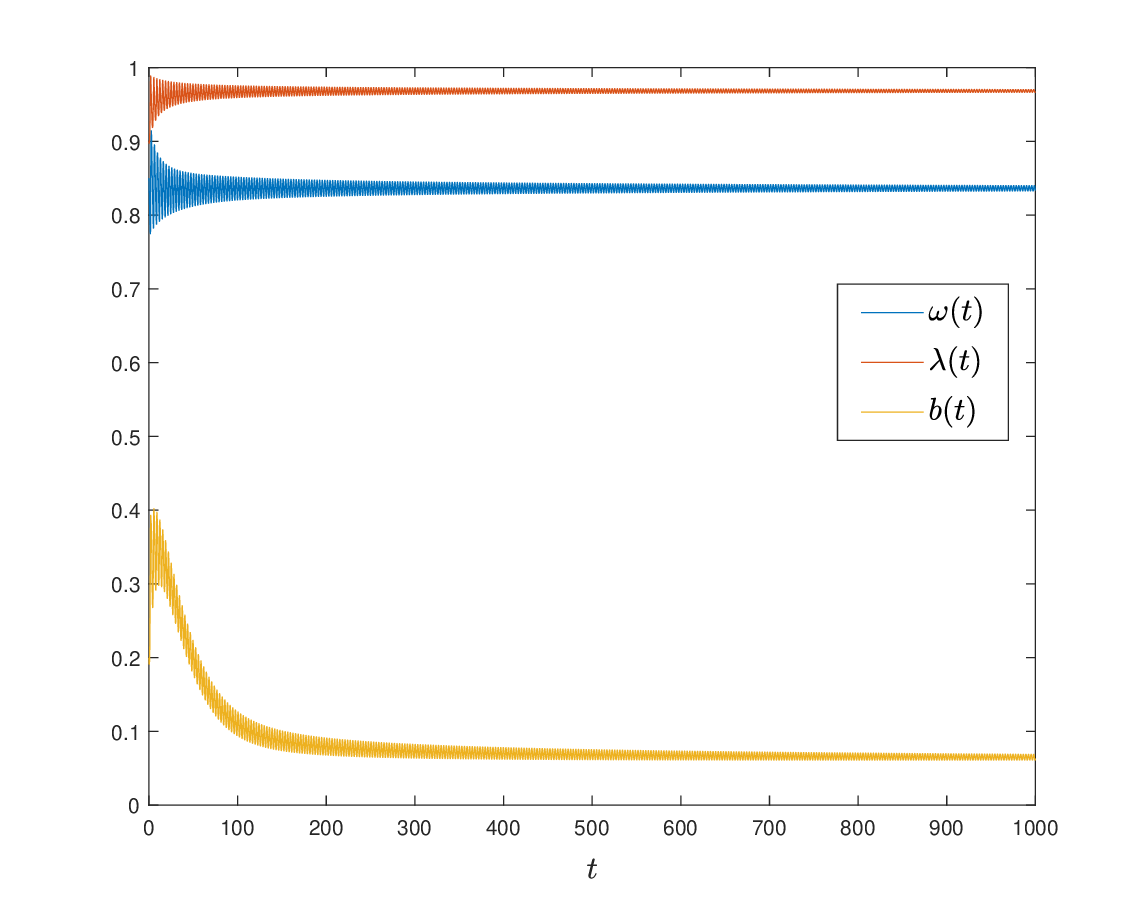}
	\caption{$\tau=\tau_0=0.82998$}
	\label{figure4}
\end{figure}
\vfill
\begin{figure}[h!]
	\centering
	\includegraphics[scale=0.60]{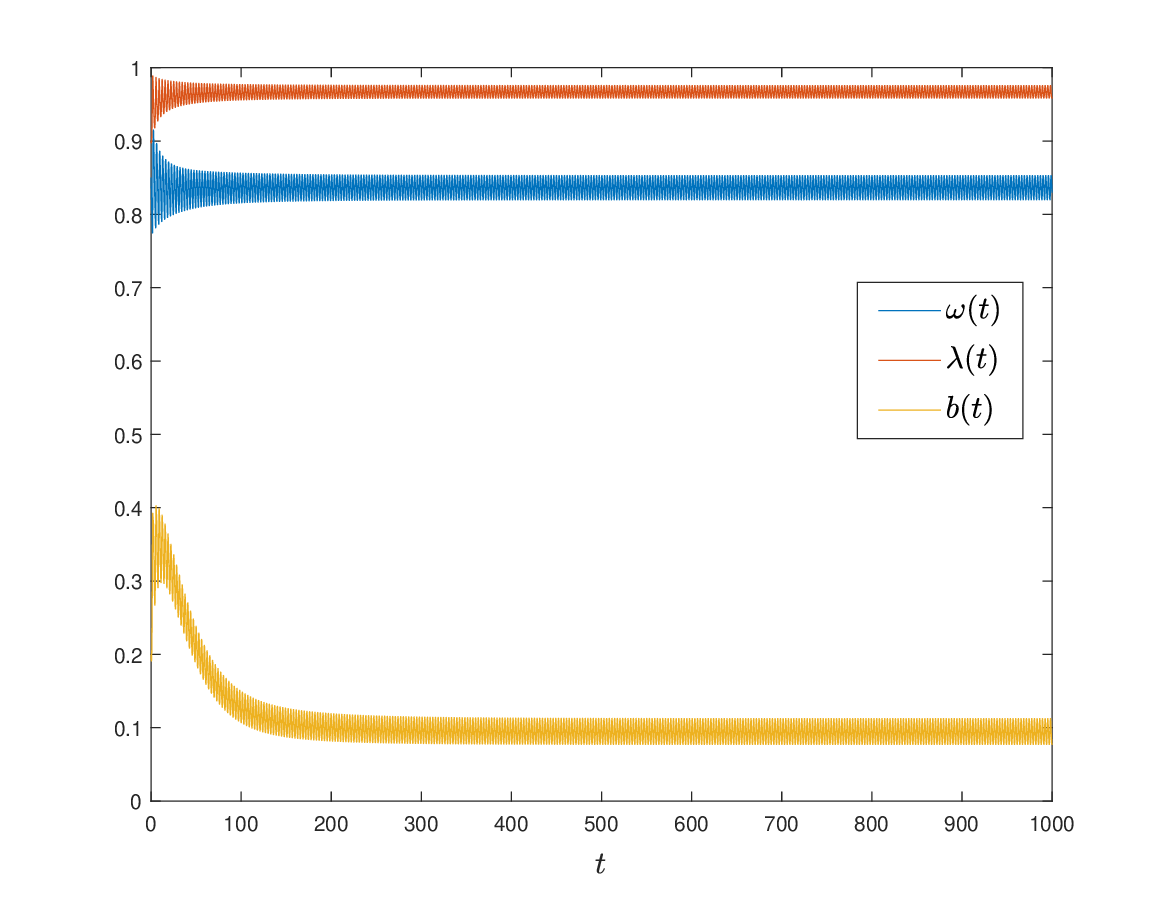}
	\caption{$\tau=0.85$}
	\label{figure4}
\end{figure}

\end{document}